\documentclass[11pt]{article}
\usepackage{amssymb,amsfonts,amsmath,latexsym,epsf,tikz,url}
\newtheorem{theorem}{Theorem}[section]

\newtheorem{corollary}[theorem]{Corollary}
\newtheorem{lemma}[theorem]{Lemma}

\newcommand{\proof}{\noindent{\bf Proof.\ }}
\newcommand{\qed}{\hfill $\square$\medskip}

\begin{document}

\title{On the some parameters related to matching of graph powers}

\author{Saeid Alikhani$^{}$\footnote{Corresponding author} \and Neda Soltani}

\date{\today}

\maketitle

\begin{center}

Department of Mathematics, Yazd University, 89195-741, Yazd, Iran\\
{\tt alikhani@yazd.ac.ir,~neda\_soltani@ymail.com}\\

\end{center}

\begin{abstract}
Let $G=(V,E)$ be a simple connected graph. A matching of $G$ is a set of disjoint edges of $G$.
For every $n, m\in\mathbb{N}$, the $n$-subdivision of $G$ is a simple graph $G^{\frac{1}{n}}$ which is constructed by replacing each edge of $G$ with a path of length $n$ and the $m$th power of $G$, denoted by $G^m$, is a graph with the same vertex set as $G$ such that two vertices are adjacent in $G^m$ if and only if their distance is at most $m$ in $G$. 
The $m^{th}$ power of the $n$-subdivision of $G$ has been introduced as a fractional power of $G$ and is denoted by 
$G^{\frac{m}{n}}$. In this paper, we study some parameters related to matching of  the natural and the fractional powers of some specific graphs. Also we study these parameters for power of  graphs that  are importance of in Chemistry. 

\end{abstract}

\noindent{\bf Keywords:} matching, subdivision of a graph, power of a graph.

\medskip
\noindent{\bf AMS Subj.\ Class.:} 05C70,  05C76.

\section{Introduction}
Let $G$ be a simple graph with vertex set $V(G)$ and edge set $E(G)$. 
The distance between every two vertices $u$ and $v$ of the graph $G$, is defined
as the length of a minimum path connecting $u$ and $v$ and is denoted by $d(u, v)$. For a vertex $u$,  
$ecc(u)=max\{d(u, x) : x\in V(G)\}$ and is called the eccentricity of $u$. The diameter $d$ of $G$ is the maximum eccentricity among vertices of $G$.

A matching $M$ in a graph $G$ is a collection of edges of $G$ such that no two edges from $M$ share a vertex. The cardinality of $M$ is called the size of the matching.
A matching $M$ is a maximum matching if there is no matching in $G$ with greater size. The cardinality of any maximum matching in $G$ is denoted by ${\alpha}^\prime(G)$ and called the matching number of $G$.
Since each vertex can be incident to at most one edge of a matching, it follows that no graph on $n$ vertices can have matching number greater than $\lfloor\frac{n}{2}\rfloor$. If each vertex of $G$ is incident with an edge of $M$, the matching $M$ is called perfect. So the number of vertices of a graph $G$ admitting a perfect matching is even, but the opposite is generally not true. Perfect matchings are obviously also maximum matchings. The study of perfect matchings, also known as Kekul\'e structures has a long history in both mathematical and chemical literature. For more details on perfect matching, we refer the reader to see \cite{Lov}. 

A matching $M$ is maximal if it cannot be extended to a larger matching in $G$. Obviously, every maximum matching is also maximal, but the opposite is generally not true. The cardinality of any smallest maximal matching in $G$ is the saturation number of $G$ and is denoted by $s(G)$. 
It is easy to see that the saturation number of a graph $G$ is at least one half of the matching number of $G$, i.e., $s(G)\geq \frac{\alpha^{\prime}(G)}{2}$ (\cite{Dos}). We recall that a set of vertices $I$ is independent if no two vertices from $I$ are adjacent. Clearly, the set of vertices that is not covered by a maximal matching is independent (\cite{Edm}). This observation provides an obvious lower bound on saturation number of the graph $G$, i.e. $s(G)\geq \frac{(|V(G)|-|I|)}{2}$ (\cite{Ves}).

For every positive integer $k$, the $k$-power of $G$ is defined on the $V(G)$ by adding edges joining any two distinct vertices $x$ and $y$ with distance at most $k$ in $G$ and is denoted by $G^k$ (\cite{{Agn},{Kra}}). In other words, $E(G^k)=\{xy: 1\leq d_G(x, y)\leq k\}$. 
The $k$-subdivision of $G$, denoted by $G^{\frac{1}{k}}$, is constructed by replacing each edge $xy$ of $G$ with a path of length $k$, say $P_{xy}$. These paths are called superedges and any new vertex is called an internal vertex or briefly 
{\it $i$-vertex} and is denoted by $(xy)_l$, if it belongs to the superedge $P_{xy}$ and has distance $l$ from the vertex 
$x$ where $l\in\{1, 2, ..., k-1\}$. Note that $(xy)_l=(yx)_{k-l}$. Also any vertex $x$ of $G^{\frac{1}{k}}$ is a terminal vertex or brifely {\it $t$-vertex}. It can be easily verified that for $k=1$, we have $G^{\frac{1}{1}}=G^1=G$ and if the graph $G$ has $p$ vertices and $q$ edges, the graph $G^{\frac{1}{k}}$ has $p+(k-1)q$ vertices and $kq$ edges.

For $m, n\in \Bbb{N}$, the fractional power of $G$ is the $m$-power of the $n$-subdivision of $G$ and is denoted by $G^{\frac{m}{n}}$. In other words, $G^{\frac{m}{n}}={(G^{\frac{1}{n}})}^m$ (\cite{Ira}).
Note that the graphs ${(G^{\frac{1}{n}})}^m$ and ${(G^m)}^{\frac{1}{n}}$ are different graphs. The fractional power of a graph has introduced by Iradmusa in \cite{Ira}. The following lemma follows from the definition of the power of graphs.

\begin{lemma}{\rm\cite{Dis}}
Let $G$ be a connected graph of order $n$ and diameter $d$. Then
\begin{enumerate}
\item[(i)] For every natural number $t \geq d$, $G^t=K_n$ and so $\alpha^{\prime}(G^t)=\lfloor\frac{n}{2}\rfloor$.
\item[(ii)] (Theorem 1 in \cite{Hob}) Let $k=ml$, where $m, l\in \Bbb{N}$. Then $G^k={(G^m)}^l$.
\item[(iii)] (Lemma 2.1 in \cite{An}) Let $x$ and $y$ be two vertices of $G$. Then 
$d_{G^k}(x, y)=\left\lceil\frac{d_G(x, y)}{k}\right\rceil$.
\end{enumerate}
\end{lemma}

In the next section, we consider the matching number of the natural powers and the fractional powers of some specific graphs. Section 3 investigates the saturation number of powers of some certain graphs that are of importance in chemistry. 

\section{Matching number of power of specific graphs}

In this section, we compute the matching number of the natural  and the fractional powers of some certain graphs such as paths, cycles, friendship graphs and complete bipartite graphs. We also present a lower and upper bound for the matching number of the fractional powers of a graph $G$.

\begin{theorem}\label{path}
Let $P_k$ and $C_k$ be a path and a cycle of order $k$, respectively. Then we have 
\begin{enumerate}
\item[(i)] For every $m\in \Bbb{N}$, ${\alpha}^\prime({P}_{k}^{m})={\alpha}^\prime({C}_{k}^{m})=\left\lfloor{\frac{k}{2}}\right\rfloor$.
\item[(ii)] For every $n\in \Bbb{N}$, ${\alpha}^\prime({P}_{k}^{\frac{1}{n}})=\left\lfloor{\frac{nk-n+1}{2}}\right\rfloor$.
\item[(iii)] For every $n\in \Bbb{N}$, ${\alpha}^\prime({C}_{k}^{\frac{1}{n}})=\left\lfloor{\frac{nk}{2}}\right\rfloor$.
\end{enumerate}
\end{theorem}

\proof
\begin{enumerate}
\item[(i)] Clearly ${\alpha}^\prime(P_k)={\alpha}^\prime(C_k)=\left\lfloor{\frac{k}{2}}\right\rfloor$. Thus we have the result by definition. 
\item[(ii)] It can be easily verified that the graph ${P}_{k}^{\frac{1}{n}}$ is a path with $n(k-1)+1$ vertices. Therefore 
$${\alpha}^\prime({P}_{k}^{\frac{1}{n}})={\alpha}^\prime(P_{n(k-1)+1})=\left\lfloor{\frac{nk-n+1}{2}}\right\rfloor.$$
\item[(iii)] Since  ${\alpha}^\prime({C}_{k}^{\frac{1}{n}})={\alpha}^\prime(C_{nk})$,  so the result follows. \qed
\end{enumerate}

The following corollary follows from the definition of the fractional power of graphs and Theorem \ref{path}.

\begin{corollary}
If $P_k$ and $C_k$ are a path and a cycle of order $k$, respectively, then for every $m, n\in \Bbb{N}$,  
\begin{enumerate}
\item[(i)] ${\alpha}^\prime({P}_{k}^{\frac{m}{n}})={\alpha}^\prime({P}_{k}^{\frac{1}{n}})$,
\item[(ii)] ${\alpha}^\prime({C}_{k}^{\frac{m}{n}})={\alpha}^\prime({C}_{k}^{\frac{1}{n}})$.
\end{enumerate}
\end{corollary}

\begin{figure}
	\begin{center}
	\includegraphics[width=10cm,height=2.3cm]{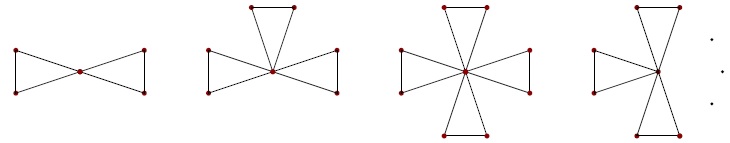}
	\caption{\small Friendship graphs $F_2, F_3, F_4$ and $F_k$, respectively.}
	\label{Fk}
\end{center}
\end{figure}

\medskip
Here we consider the matching number of friendship graphs. The friendship (or Dutch-Windmill) graph $F_k$ is a graph can be constructed by joining $k$ copies of the cycle graph $C_3$ with a common vertex. The Friendship Theorem of Paul Erd\"{o}s and et. al (\cite{Erd}) states that graphs with the property that every two vertices are connected by a path of length $2$ and which does not contain any cycle of length $4$ are friendship graphs.
Some examples of friendship graphs are shown in Figure \ref{Fk}.

\newpage
\begin{theorem}
If $F_k$ is a friendship graph, then 
 \begin{enumerate}
\item[(i)] For every $m\in \Bbb{N}$, $\alpha^{\prime}(F_k^m)=k$.
\item[(ii)] For every $n\in \Bbb{N}$, $\alpha^{\prime}(F_k^{\frac{1}{n}})=\left\{
\begin{array}{lcl}
k\big(\frac{3n-1}{2}\big), & \quad\mbox{if $n$ is odd}\\[15pt]
k\big(\frac{3n-2}{2}\big)+1, & \quad\mbox{if $n$ is even}.
\end{array}
\right.$
\item[(iii)] For every $m, n\in \Bbb{N}$, $\alpha^{\prime}(F_k^{\frac{m}{n}})=\left\lfloor\frac{k(3n-1)+1}{2}\right\rfloor$.
\end{enumerate}
\end{theorem}

\begin{figure}
	\begin{center}
	\includegraphics[width=2.7cm,height=2.7cm]{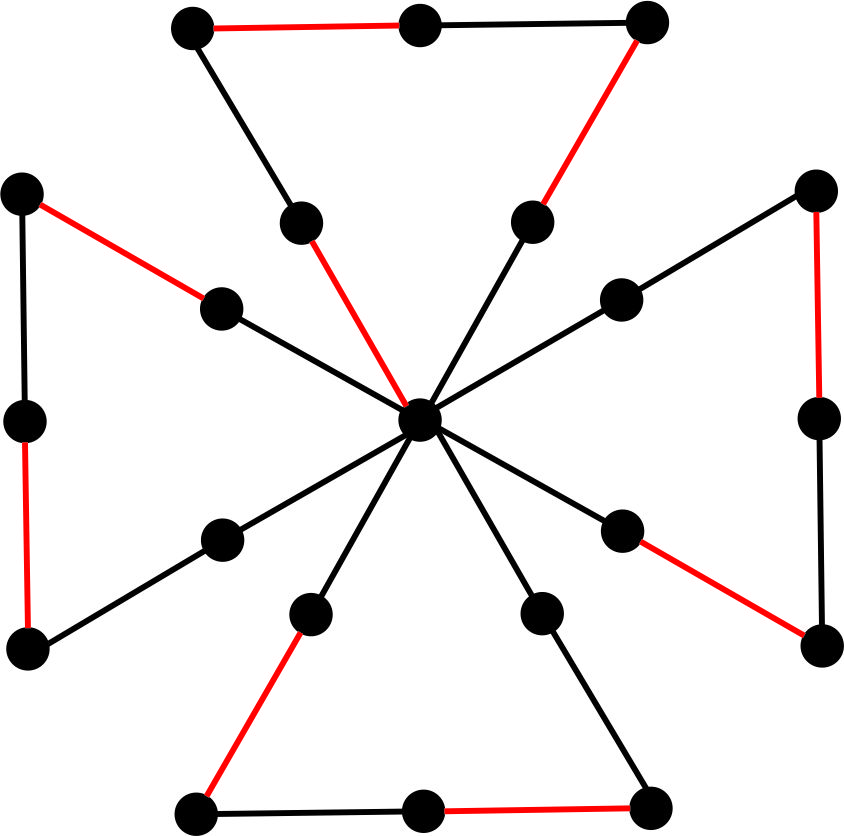}
	\caption{\small A maximum matching of the graph $F_4^{\frac{1}{2}}$.}
	\label{F4}
\end{center}
\end{figure}

\proof
 \begin{enumerate}
\item[(i)] Since for each maximum matching $M$ of the graph $F_k$, there exists only one $M$-unsaturated vertex, so for every $m\geq 2$ we have $\alpha^{\prime}(F_k^m)=\alpha^{\prime}(F_k)=k$.
\item[(ii)] Let $M$ be a maximum matching of $F_k^{\frac{1}{n}}$. Clearly, each cycle in $F_k^{\frac{1}{n}}$ has $3n$ vertices. If $n$ is odd, then the common vertex is $M$-unsaturated and so there exist $\frac{3n-1}{2}$ edges of each cycle in $M$ and this implies the result. 

Now suppose that $n$ is even. Then the number of $M$-unsaturated vertices is equal to $k-1$ (see Figure \ref{F4}) and so the result follows.
\item[(iii)] It can be easily verified that the graph $F_k^{\frac{m}{n}}$ is of order $k(3n-1)+1$. For even $n$ and odd $k$, every vertex of the graph is incident with an edge of the maximum matching of the graph and for odd $n$ or even $n, k$, there exists one vertex which is not covered by the maximum matching of $F_k^{\frac{m}{n}}$. So we have the result. 
\end{enumerate}
\qed

\begin{theorem}\label{bipart}
If $K_{m, n}$ is a complete bipartite graph and $k\in \Bbb{N}$, then
\begin{enumerate}
\item[(i)] For every $k\geq 3$, ${\alpha}^\prime(K_{m, n}^k)={\alpha}^\prime(K_{m, n}^2)=\left\lfloor{\frac{m+n}{2}}\right\rfloor$.
\item[(ii)] ${\alpha}^\prime(K_{m, n}^{\frac{1}{k}})=\left\{
\begin{array}{lcl}
(n+m)\big(\frac{k}{2}\big)+(nm-n-m)\left\lfloor{\frac{k-1}{2}}\right\rfloor, & \quad\mbox{if $k$ is even}\\[15pt]
m(\frac{k+1}{2})+m(m-1)(\frac{k-1}{2}), & \quad\mbox{if $k$ is odd and $n=m$}\\[15pt]
m(\frac{k+1}{2})+m(m-1)(\frac{k-1}{2})+m(n-m)\left\lfloor{\frac{k}{2}}\right\rfloor, & \quad\mbox{if $k$ is odd and $n>m$}.
\end{array}
\right.$
\end{enumerate}
\end{theorem}

\proof
\begin{enumerate}
\item[(i)] Since the diameter of graph $K_{m, n}$ is equal to $2$, so for every $k\geq 3$ we have 
$${\alpha}^\prime(K_{m, n}^k)={\alpha}^\prime(K_{m, n}^2).$$  
In addition, by the definition of the $k$-power of graphs, $K_{m, n}^2=K_{m+n}$ where $K_{m+n}$ is 
a complete graph with $m+n$ vertices and ${\alpha}^\prime(K_{m+n})=\left\lfloor{\frac{m+n}{2}}\right\rfloor$.
\item[(ii)] Let $M$ be a maximum matching of $K_{m, n}^{\frac{1}{k}}$. Then ${\alpha}^\prime(K_{m, n}^{\frac{1}{k}})=|M|$. If $k$ is even, then the $t$-vertices of $K_{m, n}^{\frac{1}{k}}$ can be $M$-saturated by $m+n$ edges. Note that these edges shall belong to different superedges. Now we have $m+n$ superedges which three their vertices are $M$-saturated and the other their vertices can be $M$-saturated by $\alpha^{\prime}(P_{k-2})$ edges. In addition, there exist $m(m-2)+(n-m)(m-1)$ superedges which only their end-points are $M$-saturated. So we can saturate their other vertices by $\alpha^{\prime}(P_{k-1})$ edges. Thus we have 
$${\alpha}^\prime(K_{m, n}^{\frac{1}{k}})=|M|=(n+m)\big(1+\alpha^{\prime}(P_{k-2})\big)+(nm-n-m)\alpha^{\prime}(P_{k-1})$$ and so the result follows.

Now suppose that $k$ is odd and $n=m$. In this case, we can choose $\alpha^{\prime}(P_{k+1})$ edges from $m$ superedges which are not adjacent. For the rest of $m(m-1)$ superedges, the number of edges we can put in $M$ is equal to 
$\alpha^{\prime}(P_{k-1})$. So we have $${\alpha}^\prime(K_{m, n}^{\frac{1}{k}})=|M|=m\alpha^{\prime}(P_{k+1})+m(m-1)\alpha^{\prime}(P_{k-1}).$$ and when $m<n$, we should add $(n-m)m\alpha^{\prime}(P_{k})$ edges from 
$(n-m)m$ superedges which one of their end-points is $M$-saturated. This implies the result. \qed
\end{enumerate}

\begin{corollary}
If $K_{m, n}$ is a complete bipartite graph, $k\in\Bbb{N}$ and the set $M$ is a maximum matching of $K_{m, n}^{\frac{1}{k}}$, then for every positive integer $t\geq 2$ we have
$$\alpha^{\prime}(K_{m, n}^{\frac{t}{k}})=\alpha^{\prime}(K_{m, n}^{\frac{1}{k}})+\left\lfloor{\frac{l}{2}}\right\rfloor$$ where $l$ is the number of $M$-unsaturated vertices of $K_{m, n}^{\frac{1}{k}}$.
\end{corollary}

\proof 
Suppose that $k$ is odd, $P_{xy}=xt_1t_2...t_{k-1}y$ is a superedge and $M$ is a maximum matching of $K_{m, n}^{\frac{1}{k}}$. If $xt_1, t_{k-1}y\in M$ or $xt_1, t_{k-1}y\not\in M$, then $P_{xy}$ does not have any $M$-unsaturated vertex. Otherwise, a superedge with $k-2$ vertices remains and since $k-2$ is odd, so we have one $M$-unsaturated vertex in $P_{xy}$. By the proof of Theorem \ref{bipart}, the number of these superedges is equal to $|n-m|$. Note that to have a maximum matching of $K_{m, n}^{\frac{t}{k}}$, the distance between $M$-unsaturated vertices in $K_{m, n}^{\frac{1}{k}}$ should be two, pairwise. This implies the result. 

Now assume that $k$ is even. By the proof of Theorem \ref{bipart}, there exist $nm-n-m$ superedges which their end-points are $M$-saturated. For the rest of $k-1$ vertices, we have one $M$-unsaturated vertex in each superedge. Therefore the number of $M$-unsaturated vertices in $K_{m, n}^{\frac{1}{k}}$ is equal to 
$nm-n-m$. If the distance between $M$-unsaturated vertices in $K_{m, n}^{\frac{1}{k}}$ is two, pairwise, then we have a maximum matching of $K_{m, n}^{\frac{t}{k}}$ and the result follows.
\qed

\medskip
\medskip
For any graph $G$ of order $n$, we have ${\alpha}^\prime(G)\leq {\alpha}^\prime(G^2)\leq {\alpha}^\prime(G^3)\leq\cdots\leq {\alpha}^\prime(G^d)=\left\lfloor\frac{n}{2}\right\rfloor$ where $d$ is the diameter of $G$. Because $G^k$ is an spanning subgraph of $G^{k+1}$ and so every maximum matching of $G^k$ is a matching for $G^{k+1}$. 
Also it can be easily verified that by replacing any graph by its subdivisions, the matching number of graph increases. Thus for every $k\in \Bbb{N}$, we have ${\alpha}^\prime(G)\leq {\alpha}^\prime(G^{\frac{1}{2}})\leq {\alpha}^\prime(G^{\frac{1}{3}})\leq\cdots\leq {\alpha}^\prime(G^{\frac{1}{k}})$. 
The following theorem gives a lower and upper bound for the matching number of the fractional powers of graph $G$.

\begin{theorem}
Let $G$ be a simple graph, $m, n\in \Bbb{N}$ and the set $M$ be a maximum matching of $G^{\frac{1}{n}}$. Then 
$$\alpha^{\prime}(G^{\frac{1}{n}})\leq \alpha^{\prime}(G^{\frac{m}{n}})\leq \alpha^{\prime}(G^{\frac{1}{n}})+\left\lfloor{\frac{l}{2}}\right\rfloor$$ 
where $l$ is the number of $M$-unsaturated vertices of $G^{\frac{1}{n}}$.
\end{theorem}

\proof
By the definition of the fractional power of $G$, $G^{\frac{m}{n}}={(G^{\frac{1}{n}})}^m$. The graph  $G^{\frac{1}{n}}$ is a spanning subgraph of $G^{\frac{m}{n}}$ and so every maximum matching of $G^{\frac{1}{n}}$ is a matching for $G^{\frac{m}{n}}$. Therefore $\alpha^{\prime}(G^{\frac{1}{n}})\leq \alpha^{\prime}(G^{\frac{m}{n}})$.

Now let $v_1, v_2, ..., v_l$ be  the $M$-unsaturated vertices of $G^{\frac{1}{n}}$. If for every $1\leq i\leq l$, there exists unique $1\leq j\leq l$ such that $i\neq j$ and $d(v_i, v_j)\leq m$, then we can add $\left\lfloor{\frac{l}{2}}\right\rfloor$ edges to $M$ and so $$\alpha^{\prime}(G^{\frac{m}{n}})=\alpha^{\prime}(G^{\frac{1}{n}})+\left\lfloor{\frac{l}{2}}\right\rfloor.
$$ Otherwise, the number of edges we can add to $M$ is less than $\left\lfloor{\frac{l}{2}}\right\rfloor$ and this implies the result.
\qed

\section{Saturation number of power of specific graphs}

This section investigates the saturation number of the natural and the fractional powers of some certain graphs such as paths, cycles and chain triangular cactuses. First we study the saturation number of the powers of paths and cycles.

\begin{lemma}\label{spk}
Let $P_k$ be a path with $k$ vertices, $m\in \Bbb{N}$, the set $M$ has the smallest cardinality over all maximal matchings of $P_k^m$ and $l_p$ be the number of $M$-unsaturated vertices of $P_k^m$. If $q=\left\lfloor{\frac{k}{m+1}}\right\rfloor$ and $r=k-(m+1)\left\lfloor{\frac{k}{m+1}}\right\rfloor$, then 
 
\begin{enumerate}
\item[(i)] For even $m, r$  \textbf{or} odd $m, r, q$ \textbf{or} odd $m$ and even $r, q$, we have 
$l_p=\left\lfloor{\frac{k}{m+1}}\right\rfloor$.
\item[(ii)] For even $m$ and odd $r$ \textbf{or} odd $m, r$ and even $q$ \textbf{or} odd $m, q$ and even $r$, we have $l_p=\left\lfloor{\frac{k}{m+1}}\right\rfloor+1$.
\end{enumerate}
 
\end{lemma}

\medskip
\proof
Without loss of generality, we assume that the most of edges of $M$ are the 
edges of  the path $P_k$. We should divide the path $P_k$ into $\left\lfloor{\frac{k}{m+1}}\right\rfloor$ paths of order $m+1$ which their first vertices are  $M$-unsaturated. For the rest of $k-(m+1)\left\lfloor{\frac{k}{m+1}}\right\rfloor$ vertices, we have two cases:

\medskip
Case 1) If $m$ is even, then all the vertices of the paths of order $\left\lfloor{\frac{k}{m+1}}\right\rfloor$ are $M$-saturated, except the first vertices and so the number of $M$-unsaturated vertices of $P_k^m$ depends on the path of order $k-(m+1)\left\lfloor{\frac{k}{m+1}}\right\rfloor$. If the order of this path is even, then there is no another $M$-unsaturated vertex and so $l_p=\left\lfloor{\frac{k}{m+1}}\right\rfloor$. Otherwise we have a $M$-unsaturated vertex in this path and thus $l_p=\left\lfloor{\frac{k}{m+1}}\right\rfloor+1$.

\begin{figure}[ht]
\centerline{\includegraphics[width=8cm]{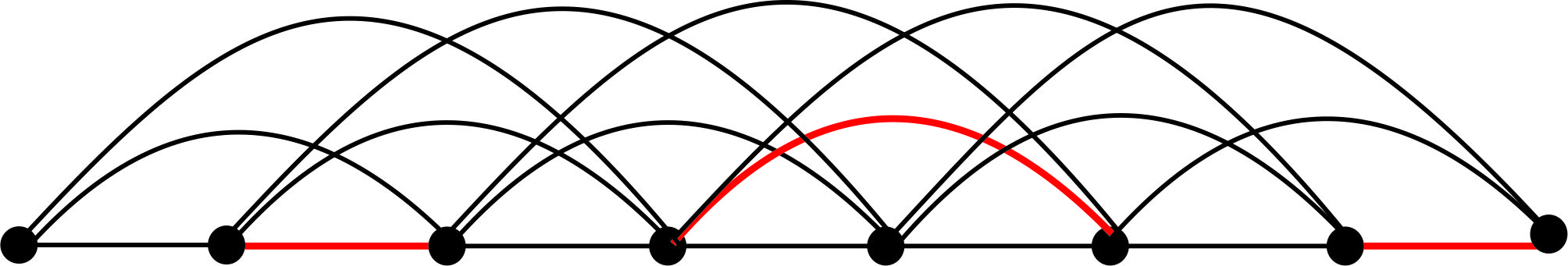}}
\caption{\label{P8}\small The smallest maximal matching of the graph $P_8^3$.}
\end{figure}

\medskip
Case 2) If $m$ is odd, then the number of vertices in the paths of order $\left\lfloor{\frac{k}{m+1}}\right\rfloor$ is even. Since the first vertex of each path is $M$-unsaturated, so we shall saturate the second $M$-unsaturated vertex of each path by the new edges in $P_k^m$ (see Figure \ref{P8}). Note that in this case, if $\left\lfloor{\frac{k}{m+1}}\right\rfloor$ is even, then  does not remain any $M$-unsaturated vertex and the amount of $l$ depends on the number of vertices of the path of order $k-(m+1)\left\lfloor{\frac{k}{m+1}}\right\rfloor$. If the number of vertices of this path is even, then $l_p=\left\lfloor{\frac{k}{m+1}}\right\rfloor$. Otherwise $l_p=\left\lfloor{\frac{k}{m+1}}\right\rfloor+1$; but if $\left\lfloor{\frac{k}{m+1}}\right\rfloor$ is odd, then we have another $M$-unsaturated vertex. Now if $k-(m+1)\left\lfloor{\frac{k}{m+1}}\right\rfloor$ is odd too, then we can saturate the $M$-unsaturated vertices by the edge between them in $P_k^m$ and so $l_p=\left\lfloor{\frac{k}{m+1}}\right\rfloor$. Otherwise $l_p=\left\lfloor{\frac{k}{m+1}}\right\rfloor+1$.
\qed

\medskip
With similar method that used in the proof of Lemma \ref{spk}, we have the following lemma which gives the number of $M$-unsaturated vertices of the $m$-power of the cycles. Note that in the cycles, the first and the last vertices are adjacent. Thus when the number of $M$-unsaturated vertices are more than $\left\lfloor{\frac{k}{m+1}}\right\rfloor$, then there exists an edge between the last 
and the first $M$-unsaturated vertices in the graph $C_k^m$ and so the number of $M$-unsaturated vertices is reduced. 

\begin{lemma}\label{sck}
Let $C_k$ be a cycle with $k$ vertices, $m\in \Bbb{N}$, the set $M$ has the smallest cardinality over all maximal matchings of $C_k^m$ and $l_c$ be the number of $M$-unsaturated vertices of $C_k^m$. If $q=\left\lfloor{\frac{k}{m+1}}\right\rfloor$ and $r=k-(m+1)\left\lfloor{\frac{k}{m+1}}\right\rfloor$, then
\begin{enumerate}
\item[(i)] For even $m, r$ \textbf{or} odd $m, r, q$ \textbf{or} odd $m$ and even $r, q$, we have $l_c=\left\lfloor{\frac{k}{m+1}}\right\rfloor$.
\item[(ii)] For even $m$ and odd $r$ \textbf{or} odd $m, r$ and even $q$ \textbf{or} odd $m, q$ and even $r$, we have $l_c=\left\lfloor{\frac{k}{m+1}}\right\rfloor-1$.
\end{enumerate}
\end{lemma}

\medskip
Now with previous lemmas, we can present the saturation number of the natural powers of the paths and the cycles. 

\begin{theorem}
Let $P_k$ and $C_k$ be a path and a cycle of order $k$, respectively and the set $M$ is a maximal matching of $G^m$ with the smallest cardinality. Then we have 

\begin{enumerate}\label{spc}
\item[(i)] For every $m\in \Bbb{N}$, $s(P_k^m)=\frac{k-l_p}{2}$.
\item[(ii)] For every $m\in \Bbb{N}$, $s(C_k^m)=\frac{k-l_c}{2}$.
\end{enumerate}
where $l_p$ and $l_c$ are the number of $M$-unsaturated vertices of $P_k^m$ and $C_k^m$, respectively.
\end{theorem}

\medskip
By the definition of the fractional power of graphs, we have $P_k^{\frac{m}{n}}=(P_k^{\frac{1}{n}})^m$ and $P_k^{\frac{1}{n}}$ is a path with $n(k-1)+1$ vertices. Also  $C_k^{\frac{1}{n}}$ is a cycle of order $nk$. 
So by replacing the path of order $n(k-1)+1$ by the path $P_k$ and the cycle of order $nk$ by the cycle $C_k$ in Theorem \ref{spc}, the saturation number of these graphs follows.

\begin{figure}[ht]
\centerline{\includegraphics[width=6cm]{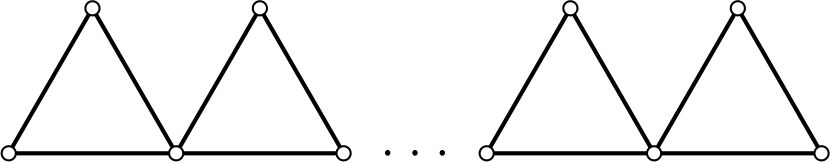}}
\caption{\label{Tk}\small A chain triangular cactus $T_k$.}
\end{figure}

\medskip
\medskip
We end this section with considering the saturation number of the natural powers and the fractional powers of some graphs with specific constructions that are of importance in chemistry. A cactus graph is a connected graph in which no edge lies in 
more than one cycle. Consequently, each block of a cactus graph is either an edge or a cycle. If all blocks of a cactus $G$ are cycles of the same size $k$, the cactus is $k$-uniform. A triangular cactus is a graph whose blocks are triangles, i.e., a $3$-uniform cactus. A vertex shared by two or more triangles is called a cut-vertex. If each triangle of a triangular cactus $G$ has at most two cut-vertices and each cut-vertex is shared by exactly two triangles, we say that $G$ is a chain triangular cactus. The number of triangles in $G$ is called the length of the chain. Obviously, all chain triangular cactus of the same length are isomorphic. Hence we denote the chain triangular cactus of length $k$ by $T_k$. An example of a chain triangular cactus is shown in Figure \ref{Tk}. Clearly, a chain triangular cactus of length $k$ has $2k+1$ vertices and $3k$ edges (\cite{Jah}). 

\begin{lemma}\label{stk}
Let $T_k$ be a chain triangular cactus of length $k$, the set $M$ has the smallest cardinality over all maximal matchings of $T_k^m$ and $l_T$ be the number of $M$-unsaturated vertices of $T_k^m$. Then we have 
$$l_T=\left\{
\begin{array}{lcl}
\left\lceil{\frac{k}{m}}\right\rceil, & \quad\mbox{if $\left\lceil{\frac{k}{m}}\right\rceil$ is odd}\\[15pt]
\left\lceil{\frac{k}{m}}\right\rceil-1, & \quad\mbox{if $\left\lceil{\frac{k}{m}}\right\rceil$ is even}.
\end{array}
\right.$$
\end{lemma} 

\proof
Since $|V(T_k^m)|=|V(T_k)|=2k+1$, so the number of $M$-unsaturated vertices of $T_k^m$ is odd. Suppose that the first $M$-unsaturated vertex of $T_k^m$ belongs to the first triangle. Then by the definition of the power of graphs, the second $M$-unsaturated vertex shall belong to $(m+1)$th triangle. Continuing this process, we have $\left\lceil{\frac{k}{m}}\right\rceil$ $M$-unsaturated vertices. Now if $\left\lceil{\frac{k}{m}}\right\rceil$ is odd, then the result follows. Otherwise, the last $M$-unsaturated vertex will be saturate and so $l_T=\left\lceil{\frac{k}{m}}\right\rceil-1$.
\qed

\begin{theorem}
If $T_k$ is a chain triangular cactus of length $k$, then
  \begin{enumerate}
           \item [(i)] For every $m\in\Bbb{N}$, $s(T_k^m)=\frac{2k-l_T+1}{2}$.
           \item[(ii)] For every $n\in\Bbb{N}$, $s(T_k^{\frac{1}{n}})=\left\{
\begin{array}{lcl}
nk, & \quad\mbox{if $n\equiv 0 ~ (mod\, 3)$}\\[15pt]
nk-\left\lfloor{\frac{k}{3}}\right\rfloor, & \quad\mbox{otherwise.}
\end{array}
\right.$ 
    \end{enumerate}
\end{theorem}

\begin{figure}[ht]\label{T3}
\centerline{\includegraphics[width=6cm]{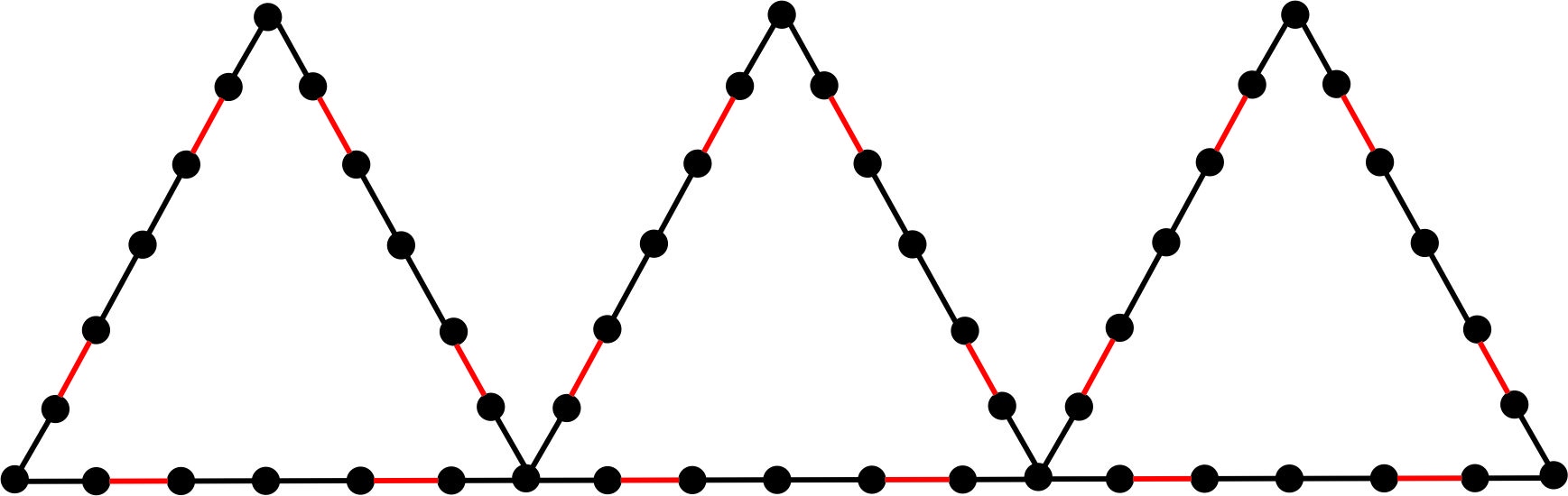}}
\caption{\label{T3}\small The smallest maximal matching of $T_3^{\frac{1}{6}}$.}
\end{figure}

\begin{figure}[ht]\label{T5}
\centerline{\includegraphics[width=8cm]{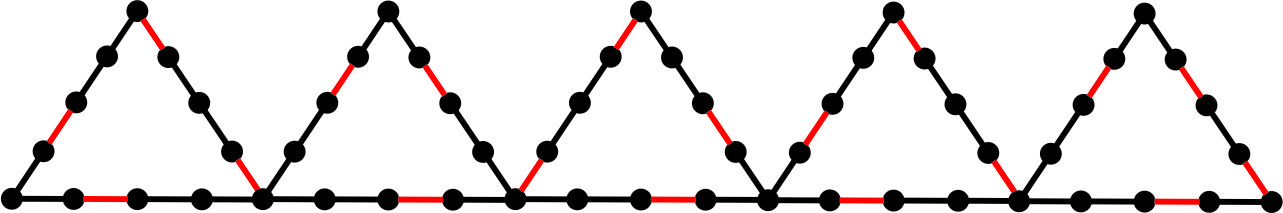}}
\caption{\label{T5}\small The smallest maximal matching of $T_5^{\frac{1}{4}}$.}
\end{figure}

\newpage
\proof
  \begin{enumerate}
           \item [(i)] It follows from Lemma \ref{stk}.
           \item[(ii)] Suppose that the set $M$ has the smallest cardinality over all maximal matchings of $T_k^{\frac{1}{n}}$. 
Clearly, each cycle in $T_k^{\frac{1}{n}}$ has $3n$ edges. If $n=3t$ where $t\in\Bbb{N}$, then we can put $n$ edges of each cycle in $M$ (see Figure \ref{T3}). So $s(T_k^{\frac{1}{n}})=|M|=nk$.

Otherwise, divide the chain $3n$-uniform cactus $T_k^{\frac{1}{n}}$ into $\left\lfloor\frac{k}{3}\right\rfloor$ chain $3n$-uniform cactuses of length $3$. Now consider the first three $3n$-uniforms. We can put $n, n-1$ and $n$ edges of these $3n$-uniforms in $M$, respectively (see Figure \ref{T5}). Continuing this process, we have $(3n-1)\left\lfloor\frac{k}{3}\right\rfloor$ edges in the maximal matching $M$. For the rest of $k-3\left\lfloor\frac{k}{3}\right\rfloor$ $3n$-uniforms, we put $n$ edges of each $3n$-uniform in $M$. Therefore $$s(T_k^{\frac{1}{n}})=|M|=(3n-1)\left\lfloor\frac{k}{3}\right\rfloor+n(k-3\left\lfloor\frac{k}{3}\right\rfloor)$$ and this implies the result.
  \end{enumerate}
\qed

\begin{theorem}\label{bound}
Let $T_k$ be a chain triangular cactus of length $k$ and $m, n\geq 2$. Then
$$nk\leq s(T_k^{\frac{m}{n}})\leq  nk+\left\lfloor\frac{nk-k+1}{2}\right\rfloor.$$
\end{theorem}

\proof
Suppose that $M$ is the smallest maximal matching of $T_k^{\frac{m}{n}}$. Then $s(T_k^{\frac{m}{n}})=|M|$. We denote the number of $M$-unsaturated vertices $T_k^{\frac{m}{n}}$ by $l$. It can be easily verified that for every constant $n$, when $m_1<m_2$, then $l_1>l_2$ and so $s(T_k^{\frac{m_1}{n}})<s(T_k^{\frac{m_2}{n}})$. Thus the lower bound of $s(T_k^{\frac{m}{n}})$ can be obtained when $m=2$.
Since $M$ has the smallest cardinality over all maximal matchings of $T_k^{\frac{m}{n}}$, the distance between a $M$-unsaturated vertex in $T_k^{\frac{2}{n}}$ with the next one should be equal to $3$. Thus there exist $3$ edges between the two $M$-unsaturated vertices in the graph and we can put the central edge in $M$. Since $T_k^{\frac{2}{n}}={(T_k^{\frac{1}{n}})}^2$ and each $3n$-uniform in $T_k^{\frac{1}{n}}$ has $3n$ edges, so there exist $n$ edges of each $3n$-uniform in $M$. Therefore $s(T_k^{\frac{2}{n}})=|M|=nk$. 

Now assume that $m$ is equal to the diameter of the graph $T_k^{\frac{m}{n}}$. Then $T_k^{\frac{m}{n}}$ is a complete graph of order $k(3n-1)+1$ and in complete graphs, every maximal matching is also maximum. Thus for every $m\leq d$, $s(T_k^{\frac{m}{n}})\leq s(T_k^{\frac{d}{n}})=\left\lfloor\frac{k(3n-1)+1}{2}\right\rfloor$ and so we have the result.
\qed

\medskip
\medskip
In the particular case, if $m+1=n$, then every $M$-unsaturated vertex is a $t$-vertices of the graph 
$T_k^{\frac{m}{n}}$ and so $s(T_k^{\frac{m}{n}})=\left\lfloor\frac{3k(n-1)+1}{2}\right\rfloor$.
Here we present the saturation number of the fractional powers of the chain triangular cactus for some special $m$ and 
$n$.

  \begin{enumerate}
           \item [(i)] $s(T_k^{\frac{3}{2}})=2k+1$.
           \item [(ii)] $s(T_k^{\frac{3}{3}})=3k+\left\lceil\frac{k}{4}\right\rceil$.
           \item [(iii)] $s(T_k^{\frac{3}{4}})=4k+\left\lceil\frac{k}{2}\right\rceil$.
           \item [(iv)] $s(T_k^{\frac{3}{5}})=\left\{
\begin{array}{lcl}
5k+3\left\lfloor\frac{k}{4}\right\rfloor+1, & \quad\mbox{if $k\equiv 1 ~ (mod\, 4)$}\\[15pt]
5k+3\left\lfloor\frac{k}{4}\right\rfloor+2, & \quad\mbox{if $k\equiv 2 ~ (mod\, 4)$}\\[15pt]
5k+3\left\lfloor\frac{k}{4}\right\rfloor+3, & \quad\mbox{if $k\equiv 3 ~ (mod\, 4)$}\\[15pt]
5k+3\left\lfloor{\frac{k-1}{4}}\right\rfloor+3, & \quad\mbox{if $k\equiv 0 ~ (mod\, 4)$.}
\end{array}
\right.$ 
           \item [(v)] $s(T_k^{\frac{4}{2}})=2k+\left\lceil\frac{k}{6}\right\rceil$.
           \item [(vi)] $s(T_k^{\frac{5}{2}})=2k+\left\lceil\frac{k+1}{4}\right\rceil$.
           \item [(vii)] $s(T_k^{\frac{5}{3}})=3k+\left\lceil\frac{k+1}{2}\right\rceil$.
           \item [(viii)] $s(T_k^{\frac{5}{4}})=4k+k$.
  \end{enumerate}

\medskip


\end{document}